\documentclass [10pt] {article}

 \usepackage [cp1251] {inputenc}
 \usepackage [english] {babel}
 \usepackage {amssymb}
 \usepackage {mathtools}
 \usepackage {amsthm}
 \usepackage [curve] {xypic}
 \usepackage {indentfirst}

 \let \$\~             
 \let \|\'             
 \def \*{\textasteriskcentered}        
 \def \0{{\rm(\*)}}
 \def \~{\tilde}       
 \def \^{\hat}
 
 \def \ge{\geqslant}
 \def \:{\colon}       
 \def \%{\times}       
 \def \[{\lvert}       
 \def \]{\rvert}
 \def \sg{\triangle}   
 \def \sx{\Delta}      
 \def \dsx{\partial\sx}
 \def \Sx{{\pmb{\boldsymbol\Delta}}}
 \def \iSx{\lefteqn{\mathring\Sx}{\phantom\Sx}}  
 \def \<{\langle}      
 \def \>{\rangle}
 \def \`{\raisebox{1pt}{$\scriptscriptstyle<$}}  
 \def \'{\raisebox{1pt}{$\scriptscriptstyle>$}}
 \def \lq{\raisebox{0.5pt}{$\scriptscriptstyle<$}}
 \def \rq{\raisebox{0.5pt}{$\scriptscriptstyle>$}}
 \def \={\sharp}
 \def \V{\vee}         

 \def \N{{\mathbb N}}
 \def \Z{{\mathbb Z}}
 \def \Q{{\mathbb Q}}
 \def \R{{\mathbb R}}
 \def \C{{\mathbb C}}

 \def \id{{\mathrm{id}}}

 \DeclareMathOperator {\Ker} {Ker}
 \let\Im\undefined
 \DeclareMathOperator {\Im} {Im}
 \DeclareMathOperator {\Hom} {Hom}
 \DeclareMathOperator {\Ext} {Ext}
 \DeclareMathOperator {\Si} {Si}

 \newcommand* {\head} [1]
 {\subsubsection* {\mathversion{bold}#1}}

 \newcommand* {\subhead} [1]
 {\addvspace\medskipamount
  \noindent {\mathversion{bold}\bf\itshape #1\/}}

 \newenvironment* {claim} [1] []
 {\begin{trivlist}\item [\hskip\labelsep {\bf #1}] \it}
 {\end{trivlist} }

 \newenvironment* {demo} [1] []
 {\begin{trivlist}\item [\hskip\labelsep {\it #1}] }
 {\end{trivlist} }

 \newenvironment* {dem}
 {\begin{trivlist}\item}
 {\end{trivlist} }

\begin {document}

 \title {\large\bf
         Straight homotopy invariants}

 \author {\normalsize\rm
          S.~S.~Podkorytov}

 \date {}

 \maketitle

 \begin {abstract} \noindent
 Let
 $X$ and $Y$ be spaces and
 $M$ be an abelian group.
 A homotopy invariant $f\:[X,Y]\to M$ is called straight if
 there exists a homomorphism $F\:L(X,Y)\to M$ such that
 $f([a])=F(\<a\>)$
 for all $a\in C(X,Y)$.
 Here
 $\<a\>\:\<X\>\to\<Y\>$ is the homomorphism induced by $a$
 between the abelian groups freely generated by $X$ and $Y$ and
 $L(X,Y)$ is a certain group of ``admissible'' homomorphisms.
 We show that all straight invariants can be expressed through
 a ``universal'' straight invariant of homological nature.
 \end {abstract}


 \head {\S~1. Introduction}


 We
 define straight homotopy invariants of maps and
 give their characterization,
 which reduces them to the classical homology theory.
 Straight invariants are a variant of the notion of homotopy
 invariants of degree at most 1
 \cite{me-3, me-1, me-2, me-4}.
 This variant has especially simple homological
 characterization.
 Homotopy invariants of finite degree are a homotopy analogue
 of Vassiliev invariants
 \cite{me-1}.


 \subhead {The group $L(X,Y)$.}
 For a set $X$,
 let $\<X\>$ be the (free) abelian group with the basis
 $X^\=\subseteq\<X\>$
 endowed with the bijection $X\to X^\=$, $x\mapsto\`x\'$.
 For sets $X$ and $Y$,
 let $L(X,Y)\subseteq\Hom(\<X\>,\<Y\>)$ be the subgroup
 generated by the homomorphisms $u$ such that
 $u(X^\=)\subseteq Y^\=\cup\{0\}$.
 (Elements of $L(X,Y)$ are the homomorphisms bounded with
 respect to the $\ell_1$-norm.)
 A map $a\:X\to Y$ induces the homomorphism
 $\<a\>\in L(X,Y)$, $\<a\>(\`x\')=\`a(x)\'$.


 \subhead {Straight homotopy invariants.}
 Let
 $X$ and $Y$ be spaces,
 $M$ be an abelian group, and
 $f\:[X,Y]\to M$ be a map (a homotopy invariant).
 The invariant $f$ is called {\it straight\/} if
 there exists a homomorphism $F\:L(X,Y)\to M$ such that
 $f([a])=F(\<a\>)$
 for all $a\in C(X,Y)$.

 (If $M$ is divisible,
 the group $L(X,Y)$ can be replaced here by
 $\Hom(\<X\>,\<Y\>)$
 because
 any homomorphism $L(X,Y)\to M$ extends to $\Hom(\<X\>,\<Y\>)$
 in this case.
 In general,
 this replacement is inadequate.
 For example,
 let $X$ and $Y$ be circles.
 Then
 the invariant ``degree'' $[X,Y]\to\Z$ is straight
 by Theorem~1.1 (or Corollary~6.8).
 At the same time,
 every homomorphism $F\:\Hom(\<X\>,\<Y\>)\to\Z$ factors through
 the restriction homomorphism
 $\Hom(\<X\>,\<Y\>)\to\Hom(\<T\>,\<Y\>)$
 for some finite set $T\subseteq X$
 \cite[\S~94]{Fuchs}.
 Thus $F$ cannot give rise to a non-constant homotopy
 invariant.)


 \subhead {The main invariant $h\:[X,Y]\to[SX,SY]$.}
 For a space $X$,
 let $SX$ be its singular chain complex.
 Let $X$ and $Y$ be spaces.
 Let $[SX,SY]$ be the group of homotopy classes of morphisms
 $SX\to SY$.
 There is a (non-naturally) split exact natural sequence
 $$
 0
 \longrightarrow
 \prod_{i\in\Z}\Ext(H_{i-1}X,H_iY)
 \longrightarrow
 [SX,SY]
 \longrightarrow
 \prod_{i\in\Z}\Hom(H_iX,H_iY)
 \longrightarrow
 0
 $$
 (``the universal coefficient theorem'',
 cf.\ \cite[Theorem~5.5.3]{Spanier}).
 For $a\in C(X,Y)$,
 let
 $Sa\:SX\to SY$ be the induced morphism and
 $[Sa]\in[SX,SY]$ be its homotopy class.
 The invariant $h\:[X,Y]\to[SX,SY]$, $[a]\mapsto[Sa]$, is
 called {\it main}.


 \subhead {The main result.}
 We call a space {\it valid\/} if
 it is homotopy equivalent to a CW-complex;
 we call it {\it finitary\/} if
 it is weakly homotopy equivalent to a compact CW-complex.

 \begin {claim} [1.1. Theorem.]
 Let
 $X$ be a finitary valid space,
 $Y$ be a valid space,
 $h\:[X,Y]\to[SX,SY]$ be the main invariant,
 $M$ be an abelian group, and
 $f\:[X,Y]\to M$ be an invariant.
 The invariant $f$ is straight
 if and only if
 there exists a homomorphism $d\:[SX,SY]\to M$ such that
 $f=d\circ h$.
 \end {claim}

 \begin {dem}
 The theorem follows from Propositions 7.3 and 12.2.
 \qed
 \end {dem}

 The theorem says that
 the main invariant is a ``universal'' straight invariant.
 A weaker and slightly complicated result is
 \cite[Theorem~II]{me-0}.
 If $M$ is divisible,
 then
 the sufficiency (``if'') follows easily from an appropriate
 form of the Dold--Thom theorem
 (see \S~7), and
 the necessity (``only if'') follows from
 \cite[Theorem~II]{me-0}
 (but any abelian group is a subgroup of a divisible one).
 The validity and finitarity hypotheses are essential,
 see \S\S\ 13, 14.

 In \S~15,
 we consider $K$-straight invariants taking values in modules
 over a commutative ring $K$
 (by definitions,
 straight = $\Z$-straight).


 \head {\S~2. Notation}


 \subhead {The question mark.}
 The expression $[?]$ denotes the map $a\mapsto[a]$ between
 sets indicated in the context.
 We similarly use $\<?\>$, etc.
 This notation is also used for functors.


 \subhead {Sets and abelian groups.}
 For a set $X$,
 let $c_X\:X\to\<X\>$ be the canonical map $x\mapsto\`x\'$.
 For $v\in\<X\>$ and $x\in X$,
 let $v/x\in\Z$ be the coefficient of $\`x\'$ in $v$.
 For an abelian group $T$,
 a map $a\:X\to T$ gives rise to the homomorphism
 $a^+\:\<X\>\to T$, $\`x\'\mapsto a(x)$.
 $T^X$ is the group of maps $X\to T$.


 \subhead {Simplicial sets.}
 For simplicial sets $U$ and $V$,
 let
 $\Si(U,V)$ be the set of simplicial maps and
 $[U,V]$ be the set of their homotopy classes
 (two simplicial maps are homotopic if
 they are connected by a sequence of homotopies).
 The functor $\<?\>$ takes simplicial sets to simplicial
 abelian groups degreewise.
 There is the canonical simplicial map $c_U\:U\to\<U\>$.
 For a simplicial abelian group $Z$,
 a simplicial map $s\:U\to Z$ gives rise to the simplicial
 homomorphism $s^+\:\<U\>\to Z$.
 For a simplicial set $T$,
 a simplicial map $s\:U\to V$ induces the maps
 $s^T_\#\:\Si(T,U)\to\Si(T,V)$,
 $s_T^\#\:\Si(V,T)\to\Si(U,T)$,
 $s^T_*\:[T,U]\to[T,V]$, and
 $s_T^*\:[V,T]\to[U,T]$.
 This notation is also used in the topological case.
 

 \head {\S~3. Induced straight invariants}


 \begin {claim} [3.1. Lemma.]
 Let
 $X$, $\~X$, $\~Y$, and $Y$ be spaces,
 $r\:X\to\~X$ and $s\:\~Y\to Y$ be continuous maps,
 $M$ be an abelian group and
 $f\:[X,Y]\to M$ be a straight invariant.
 Then the invariant $\~f\:[\~X,\~Y]\to M$,
 $\~f([\~a])=f([s\circ\~a\circ r])$, $\~a\in C(\~X,\~Y)$, is
 straight.
 \end {claim}

 \begin {demo} [Proof.]
 There is a homomorphism $F\:L(X,Y)\to M$ such that
 $f([a])=F(\<a\>)$, $a\in C(X,Y)$.
 We have the commutative diagram
 $$
 \xymatrix {
 C(\~X,\~Y)
 \ar[rr]^-{\<?\>}
 \ar[rd]^-{K}
 \ar[dd]_-{[?]} & &
 L(\~X,\~Y)
 \ar[d]_-{T}
 \ar@/^6ex/[dd]^-{\~F} \\
 &
 C(X,Y)
 \ar[r]^-{\<?\>}
 \ar[d]_-{[?]} &
 L(X,Y)
 \ar[d]_-{F} \\
 [\~X,\~Y]
 \ar[r]^-{k}
 \ar@/_3ex/[rr]_-{\~f} &
 [X,Y]
 \ar[r]^-{f} &
 M,
 }
 $$
 where
 the maps $K$ and $k$ and
 the homomorphism $T$
 are induced by the pair $(r,s)$
 (that is,
 $K(\~a)=s\circ\~a\circ r$,
 $k([\~a])=[s\circ\~a\circ r]$,
 $T(\~u)=\<s\>\circ\~u\circ\<r\>$), and
 $\~F=F\circ T$.
 Thus $\~f$ is straight.
 \qed
 \end {demo}


 \head {\S~4. The main invariant
        $h\:[\[U\],\[V\]]\to[S\[U\],S\[V\]]$}


 The geometric realization $\[Z\]$ of a simplicial abelian
 group $Z$ has a structure of an abelian group.
 $\[Z\]$ is a topological abelian group if
 $Z$ is countable;
 in general,
 it is a group of the category of compactly generated
 Hausdorff spaces.
 For a simplicial set $T$,
 $C(\[T\],\[Z\])$ and $[\[T\],\[Z\]]$ are abelian groups with
 respect to pointwise addition.
 Clearly,
 $\Si(T,Z)$ and $[T,Z]$ are also abelian groups.
 
 \begin {claim} [4.1 Lemma.]
 Let $U$ and $V$ be simplicial sets.
 Then there exists a commutative diagram
 $$
 \xymatrix {
 [U,V]
 \ar[rr]^-{(c_V)^U_*}
 \ar[dd]_-{i} & &
 [U,\<V\>]
 \ar[dd]^-{j}
 \ar[dl]^-{e} \\
 &
 [S\[U\],S\[V\]]
 \ar[dr]^-{E} &
 \\
 [\[U\],\[V\]]
 \ar[rr]^-{\[c_V\]^{\[U\]}_*}
 \ar[ur]^-{h} & &
 [\[U\],\[\<V\>\]],
 }
 $$
 where
 $i\:[s]\mapsto[\[s\]]$ (the map induced by the geometric
 realization map),
 $j$ is similar,
 $h$ is the main invariant, and
 $e$, $E$ are some isomorphisms.
 \end {claim}

 This is a version of the Dold--Thom theorem
 \cite[\S~4.K]{Hatcher}.

 \begin {demo} [Proof.]
 Let $\sg$ be the singular functor.
 For a simplicial set $T$,
 let $k_T\:T\to\sg\[T\]$ be the canonical weak equivalence.
 If $T$ is a simplicial abelian group,
 $k_T$ is a simplicial homomorphism.
 We have the commutative diagram
 $$
 \xymatrix {
 V
 \ar[rr]^-{c_V}
 \ar[dd]_-{k_V} & &
 \<V\>
 \ar[dd]^-{k_{\<V\>}}
 \ar[dl]^-{\<k_V\>} \\
 &
 \<\sg\[V\]\>
 \ar[dr]^-{m} &
 \\
 \sg\[V\]
 \ar[rr]^-{\sg\[c_V\]}
 \ar[ur]^-{c_{\sg\[V\]}} & &
 \sg\[\<V\>\],
 }
 $$
 where $m=(\sg\[c_V\])^+$.
 $k_{\<V\>}$, $\<k_V\>$, and thus $m$ are weak equivalences.
 Consider the commutative diagram
 $$
 \xymatrix {
 [U,V]
 \ar[rr]^-{(c_V)^U_*}
 \ar[dd]_-{(k_V)^U_*}
 \ar@/_8ex/[ddd]_-{i} & &
 [U,\<V\>]
 \ar[dd]_-{(k_{\<V\>})^U_*}
 \ar@/^8ex/[ddd]^-{j}
 \ar[dl]^-{\<k_V\>^U_*} \\
 &
 [U,\<\sg\[V\]\>]
 \ar[dr]^-{m^U_*} &
 \\
 [U,\sg\[V\]]
 \ar[rr]^-{(\sg\[c_V\])^U_*}
 \ar[ur]^-{(c_{\sg\[V\]})^U_*} & &
 [U,\sg\[\<V\>\]] \\
 [\[U\],\[V\]]
 \ar[rr]^-{\[c_V\]^{\[U\]}_*}
 \ar[u]_-{p} & &
 [\[U\],\[\<V\>\]],
 \ar[u]^-{q}
 }
 $$
 where
 the upper part is the result of applying the functor $[U,?]$
 to the previous diagram and
 $p$ and $q$ are the standard adjunction bijections for the
 functors $\[?\]$ and $\sg$.
 $\<k_V\>^U_*$, $m^U_*$, and $q$ are isomorphisms.

 We will find an isomorphism
 $P\:[S\[U\],S\[V\]]\to[U,\<\sg\[V\]\>]$ such that
 $P\circ h=(c_{\sg\[V\]})^U_*\circ p$.
 Then it will be enough to set
 $e=P^{-1}\circ\<k_V\>^U_*$ and
 $E=q^{-1}\circ m^U_*\circ P$.

 For a simplicial set $T$,
 let $AT$ be its chain complex,
 so that $(AT)_n=\<T_n\>$, $n\ge0$.
 Then $SX=A\sg X$
 for any space $X$.
 A simplicial map $s\:T\to\<W\>$ gives rise to the morphism
 $v\:AT\to AW$, $v_n=s_n^+$, $n\ge0$.
 This rule yields an isomorphism $d\:[T,\<W\>]\to[AT,AW]$
 (the Dold--Kan correspondence).
 We set
 $T=\sg\[U\]$ and
 $W=\sg\[V\]$.
 Consider the commutative diagram
 $$
 \xymatrix {
 [\[U\],\[V\]]
 \ar[r]^-{b}
 \ar@/^5ex/[rrr]^-{p}
 \ar[d]_-{h} &
 [\sg\[U\],\sg\[V\]]
 \ar[rr]^-{(k_U)_{\sg\[V\]}^*}
 \ar[d]_-{(c_{\sg\[V\]})^{\sg\[U\]}_*} &&
 [U,\sg\[V\]]
 \ar[d]^-{(c_{\sg\[V\]})^U_*} \\
 [A\sg\[U\],A\sg\[V\]]
 \ar@/_3ex/[rrr]_-{P} &
 [\sg\[U\],\<\sg\[V\]\>]
 \ar[l]_-{d}
 \ar[rr]^-{(k_U)_{\<\sg\[V\]\>}^*} &&
 [U,\<\sg\[V\]\>],
 }
 $$
 where
 the map $b$ is given by the functor $\sg$ and
 $P=(k_U)_{\<\sg\[V\]\>}^*\circ d^{-1}$.
 Since $(k_U)_{\<\sg\[V\]\>}^*$ is an isomorphism,
 $P$ is an isomorphism too.
 \qed
 \end {demo}


 \head {\S~5. N\"obeling--Bergman theory}


 By a {\it ring\/} we mean a (non-unital) commutative ring;
 {\it subring\/} is understood accordingly.
 The following facts follow from
 \cite[Theorem~2 and its proof]{Hill},
 cf.\ \cite[\S~97]{Fuchs}.

 \begin {claim} [5.1. Lemma.]
 Let $E$ be a torsion-free ring generated by idempotents.
 Then $E$ is a free abelian group.
 \qed
 \end {claim}

 An example:
 the ring $B(X)$ of bounded functions $X\to\Z$,
 where $X$ is an arbitrary set.

 \begin {claim} [5.2. Lemma.]
 Let
 $E$ be a torsion-free ring and
 $F\subseteq E$ be a subring,
 both generated by idempotents.
 Then the abelian group $E/F$ is free.
 \qed
 \end {claim}

 For $F=0$,
 this is Lemma~5.1.


 \head {\S~6. Maps to a space with addition}


 Let
 $X$ be a space and
 $T$ be a Hausdorff space.

 For a set $V\subseteq T$,
 we introduce the homomorphism
 $s_V\:L(X,T)\to\Z^X$, $s_V(u)(x)=I_V^+(u(\`x\'))$, $x\in X$,
 where $I_V\:T\to\Z$ is the indicator function of the set $V$.


 \subhead {The subgroup $R\subseteq L(X,T)$.}
 For $p\in X$, $q\in T$,
 let $R(p,q)\subseteq L(X,T)$ be the subgroup of homomorphisms
 $u$ such that,
 for any sufficiently small (open) neighbourhood $V$ of $q$,
 the function $s_V(u)$ is constant in some neighbourhood of
 $p$.
 Let $R\subseteq L(X,T)$ be the intersection of the subgroups
 $R(p,q)$, $p\in X$, $q\in T$.

 \begin {claim} [6.1. Lemma.]
 For $a\in C(X,T)$,
 we have $\<a\>\in R$.
 \end {claim}

 \begin {demo} [Proof.]
 Take $p\in X$, $q\in T$.
 We show that $\<a\>\in R(p,q)$.
 If $a(p)=q$,
 then,
 for any neighbourhood $V$ of $q$,
 we
 take the neighbourhood $U=a^{-1}(V)$ of $p$ and
 get $s_V(\<a\>)|_U=1$.
 Otherwise,
 choose disjoint neighbourhoods
 $W$ of $q$ and
 $W_1$ of $a(p)$.
 Consider the neighbourhood $U=a^{-1}(W_1)$ of $p$.
 For any $V\subseteq W$,
 we have $s_V(\<a\>)|_U=0$.
 \qed
 \end {demo}

 \begin {claim} [6.2. Lemma.]
 The abelian group $L(X,T)/R$ is free.
 \end {claim}

 \begin {demo} [Proof.]
 Let $O_T$ be the set of open sets in $T$.
 Consider the ring $E=B(X\times X\times O_T)$.
 For $p\in X$, $q\in T$,
 let $I(p,q)\subseteq E$ be the ideal of functions $f$ such
 that,
 for any sufficiently small neighbourhood $V$ of $q$,
 the function $X\to\Z$, $x\mapsto f(p,x,V)$, vanishes in some
 neighbourhood of $p$.
 Let $I\subseteq E$ be the intersection of the ideals $I(p,q)$,
 $p\in X$, $q\in T$.
 The ring $E/I$ is
 torsion-free and
 generated by idempotents.
 By Lemma~5.1,
 $E/I$ is a free abelian group.
 Consider the homomorphism $k\:L(X,T)\to E$,
 $k(u)(p,x,V)=s_V(u)(x)-s_V(u)(p)$, $p,x\in X$, $V\in O_T$,
 $u\in L(X,T)$.
 We have
 $k^{-1}(I(p,q))=R(p,q)$ and
 thus $k^{-1}(I)=R$.
 Therefore,
 $k$ induces a monomorphism $L(X,T)/R\to E/I$.
 It follows that
 the abelian group $L(X,T)/R$ is free.
 \qed
 \end {demo}


 \subhead {The set $Q$ and the homomorphisms $e(D,a)$.}
 Let $Q$ be the set of pairs $(D,a)$,
 where
 $D\subseteq X$ is a closed set and
 $a\in C(D,T)$.
 For $(D,a)\in Q$,
 introduce the homomorphism $e(D,a)\in L(X,T)$,
 $$
 e(D,a)(\`x\')=
 \begin{cases}
 \`a(x)\' &
 \text{if $x\in D$,} \\
 0 &
 \text{otherwise,}
 \end{cases}
 $$
 $x\in X$.

 \begin {claim} [6.3. Lemma.]
 Let
 $(D,a)\in Q$,
 $p\in X$, and
 $q\in T$.
 If $e(D,a)\notin R(p,q)$,
 then $p\in D$ and $a(p)=q$.
 \end {claim}

 \begin {demo} [Proof.]
 Put $u=e(D,a)$.
 {\it The case $p\notin D$.\/}
 Consider the neighbourhood $U=X\setminus D$ of $p$.
 We have $s_V(u)|_U=0$
 for any $V\subseteq T$.
 Thus $u\in R(p,q)$.
 {\it The case $p\in D$, $a(p)\ne q$.\/}
 Choose disjoint neighbourhoods
 $W$ of $q$ and
 $W_1$ of $a(p)$.
 There is a neighbourhood $U$ of $p$ such that
 $a(D\cap U)\subseteq W_1$.
 We have $s_V(u)|_U=0$
 for any $V\subseteq W$.
 Thus $u\in R(p,q)$.
 \qed
 \end {demo}


 \subhead {The subgroup $K\subseteq L(X,T)$.}
 Let $K\subseteq L(X,T)$ be the subgroup generated by $e(D,a)$,
 $(D,a)\in Q$.

 \begin {claim} [6.4. Lemma.]
 The abelian group $L(X,T)/K$ is free.
 \end {claim}

 \begin {demo} [Proof.]
 Consider the monomorphism $j\:L(X,T)\to B(X\times T)$,
 $j(u)(x,t)=u(\`x\')/t$.
 For $(D_i,a_i)\in Q$, $i=1,2$,
 we have $j(e(D_1,a_1))j(e(D_2,a_2))=j(e(D,a))$,
 where
 $D=\{x\in D_1\cap D_2:a_1(x)=a_2(x)\}$ and
 $a=a_1|_D=a_2|_D$.
 In particular,
 $j(e(D,a))$, $(D,a)\in Q$, are idempotents.
 Therefore,
 $j(K)$ is a subring generated by idempotents.
 By Lemma~5.2,
 the abelian group $B(X\times T)/j(K)$ is free.
 Since $j$ induces a monomorphism
 $L(X,T)/K\to B(X\times T)/j(K)$,
 the abelian group $L(X,T)/K$ is free.
 \qed
 \end {demo}

 \begin {claim} [6.5. Lemma.]
 The abelian group $L(X,T)/(K\cap R)$ is free.
 \end {claim}

 \begin {demo} [Proof.]
 The quotients in the chain
 $L(X,T)\supseteq K\supseteq K\cap R$ are free:
 $L(X,T)/K$ by Lemma~6.4, and
 $K/(K\cap R)$ as a subgroup of $L(X,T)/R$,
 which is free by Lemma~6.2.
 \qed
 \end {demo}


 \subhead {The homomorphism $G\:L(X,T)\to T^X$.}
 Let $T$ have a structure of an abelian group such that,
 \0
 for any closed set $D\subseteq X$,
 the set $C(D,T)$ becomes an abelian group with respect to
 pointwise addition\footnotemark[1].
 \footnotetext[1]
 {The condition \0 is satisfied
 if $T$ is a topological abelian group or
 if $X=\[U\]$ and $T=\[Z\]$,
 where
 $U$ is a simplicial set and
 $Z$ is a simplicial abelian group.}
 Consider the homomorphism $G\:L(X,T)\to T^X$,
 $G(u)(x)=r(u(\`x\'))$, $x\in X$, $u\in L(X,T)$,
 where $r=\id^+\:\<T\>\to T$.

 \begin {claim} [6.6. Lemma.]
 $G(K\cap R)\subseteq C(X,T)$.
 \end {claim}

 \begin {demo} [Proof.]
 Take $u\in K\cap R$.
 We show that
 $G(u)\in C(X,T)$.
 Since $u\in K$,
 we have
 \begin {align*}
 u&=\sum_{i\in I}u_i, &
 u_i&=k_ie(D_i,a_i),
 \end {align*}
 where
 $I$ is a finite set,
 $k_i\in\Z$, and
 $(D_i,a_i)\in Q$.
 For $J\subseteq I$,
 put
 \begin {align*}
 u_J&=\sum_{i\in J}u_i, &
 D_J&=\bigcap_{i\in J}D_i\subseteq X
 \end {align*}
 (so $D_\varnothing=X$)
 and
 \begin {align*}
 b_J&=\sum_{i\in J}k_ia_i|_{D_J}\in C(D_J,T), &
 k_J&=\sum_{i\in J}k_i.
 \end {align*}

 Take $p\in X$.
 We verify that
 $G(u)$ is continuous at $p$.
 Put $N=\{i\in I:p\notin D_i\}$.
 For $q\in T$,
 put $I(q)=\{i\in I:p\in D_i,\ a_i(p)=q\}$.
 We have
 $$
 u=u_N+\sum_{q\in T}u_{I(q)}
 $$
 (almost all summands are zero).
 Clearly,
 $G(u_N)$ vanishes in some neighbourhood of $p$.
 Take $q\in T$.
 It suffices to show that
 $G(u_{I(q)})$ is continuous at $p$.
 Put $t_0=G(u_{I(q)})\in T$.
 We have $t_0=k_{I(q)}q$.
 Let $W$ be a neighbourhood of $t_0$.
 We seek a neighbourhood $U$ of $p$ such that
 $G(u_{I(q)})(U)\subseteq W$.

 Put $E=\{J\subseteq I(q):k_J=k_{I(q)}\}$.
 For $J\in E$,
 we have
 $p\in D_J$ and
 $b_J(p)=t_0$.
 There is a neighbourhood $U_1$ of $p$ such that
 $b_J(D_J\cap U_1)\subseteq W$
 for all $J\in E$.

 By Lemma~6.3,
 $u_i\in R(p,q)$
 for $i\in I\setminus I(q)$.
 Since $u\in R(p,q)$,
 we have $u_{I(q)}\in R(p,q)$.
 Therefore,
 there is a neighbourhood $V\subseteq T$ of $q$ such that
 the function $s_V(u_{I(q)})$ is constant
 in some neighbourhood $U_2$ of $p$.

 There is a neighbourhood $U_3$ of $p$ such that
 $a_i(D_i\cap U_3)\subseteq V$
 for all $i\in I(q)$.
 For $x\in X$,
 put $J(x)=\{i\in I(q):x\in D_i\}$.
 For $x\in U_2\cap U_3$,
 we have $k_{J(x)}=s_V(u_{I(q)})(x)=s_V(u_{I(q)})(p)=k_{I(q)}$,
 i.~e.\ $J(x)\in E$.

 Set $U=U_1\cap U_2\cap U_3$.
 Take $x\in U$.
 We have $G(u_{I(q)})(x)=b_{J(x)}(x)\in W$
 because $J(x)\in E$.
 \qed
 \end {demo}

 \begin {claim} [6.7. Lemma.]
 There exists a homomorphism $g\:L(X,T)\to C(X,T)$ such that
 $g(\<a\>)=a$
 for all $a\in C(X,T)$.
 \end {claim}

 \begin {demo} [Proof.]
 We have $G(\<a\>)=a$
 for all $a\in T^X$.
 Since
 $G(K\cap R)\subseteq C(X,T)$
 (by Lemma~6.6) and
 the abelian group $L(X,T)/(K\cap R)$ is free
 (by Lemma~6.5),
 there is a homomorphism $g\:L(X,T)\to C(X,T)$ such that
 $g(u)=G(u)$
 for $u\in K\cap R$.
 For $a\in C(X,T)$,
 we have
 $\<a\>\in K$
 (because $\<a\>=e(X,a)$)
 and
 $\<a\>\in R$
 (by Lemma~6.1).
 We get $g(\<a\>)=G(\<a\>)=a$.
 \qed
 \end {demo}

 \begin {claim} [6.8. Corollary.]
 Suppose that
 \0
 $[X,T]$ is an abelian group with respect to pointwise
 addition\footnotemark[2].
 \footnotetext[2]
 {See footnote~1.}
 Then the invariant $\id\:[X,T]\to[X,T]$ is straight.
 \end {claim}

 \begin {demo} [Proof.]
 By Lemma~6.7,
 there is a homomorphism $g\:L(X,T)\to C(X,T)$ such that
 $g(\<a\>)=a$
 for all $a\in C(X,T)$.
 Consider the homomorphism $F\:L(X,T)\to[X,T]$,
 $u\mapsto[g(u)]$.
 For $a\in C(X,T)$,
 we have $[a]=[g(\<a\>)]=F(\<a\>)$.
 \qed
 \end {demo}


 \head {\S~7. Sufficiency in Theorem~1.1}


 The proof of sufficiency in Theorem~1.1 relies on
 Corollary~6.8.
 If the group $M$ is divisible,
 it is easy to use Lemma~7.1 instead
 (then the stuff of \S\S~5, 6 is needless).

 \begin {claim} [7.1. Lemma
                 \rm (cf.\ {\cite[Lemma~1.2]{me-3}}).]
 Let
 $X$ and $T$ be spaces and
 $T$ have a structure of an abelian group such that
 \0
 the sets $C(X,T)$ and $[X,T]$ become abelian groups with
 respect to pointwise addition\footnotemark[3].
 \footnotetext[3]
 {See footnote~1.}
 Let
 $M$ be a divisible abelian group and
 $f\:[X,T]\to M$ be a homomorphism.
 Then $f$ is a straight invariant.
 \end {claim}

 \begin {demo} [Proof.]
 Consider the homomorphism $G\:L(X,T)\to T^X$,
 $G(u)(x)=r(u(\`x\'))$, $x\in X$, $u\in L(X,T)$,
 where $r=\id^+\:\<T\>\to T$.
 Let $D\subseteq L(X,T)$ be the subgroup generated by the
 homomorphisms $\<a\>$, $a\in C(X,T)$.
 Clearly,
 $G(\<a\>)=a$
 for $a\in C(X,T)$.
 Therefore,
 $G(D)\subseteq C(X,T)$.
 Consider the homomorphism $F_0\:D\to M$, $u\mapsto f([G(u)])$.
 Since $M$ is divisible,
 there is a homomorphism $F\:L(X,T)\to M$ such that
 $F|_D=F_0$.
 For $a\in C(X,T)$,
 we have $f([a])=f([G(\<a\>)])=F_0(\<a\>)=F(\<a\>)$.
 \qed
 \end {demo}

 \begin {claim} [7.2. Claim.]
 Let $U$ and $V$ be simplicial sets.
 Then the main invariant $h\:[\[U\],\[V\]]\to[S\[U\],S\[V\]]$
 is straight.
 \end {claim}

 \begin {demo} [Proof.]
 Consider the commutative diagram
 $$
 \xymatrix {
 [\[U\],\[V\]]
 \ar[r]^-{h}
 \ar[dr]_-{\[c_V\]^{\[U\]}_*}
 &
 [S\[U\],S\[V\]]
 \ar[d]^-{E} \\
 &
 [\[U\],\[\<V\>\]],
 }
 $$
 where $E$ is the isomorphism from Lemma~4.1.
 By Corollary~6.8,
 the invariant $\id\:[\[U\],\[\<V\>\]]\to[\[U\],\[\<V\>\]]$ is
 straight.
 Therefore,
 by Lemma~3.1,
 the invariant $\[c_V\]^{\[U\]}_*$ is straight.
 Since $E$ is an isomorphism,
 $h$ is also straight.
 \qed
 \end {demo}

 \begin {claim} [7.3. Proposition.]
 Let
 $X$ be a space and
 $Y$ be a valid space.
 Then the main invariant $h\:[X,Y]\to[SX,SY]$ is straight.
 \end {claim}

 \begin {demo} [Proof.]
 There are homology equivalences
 $r\:\[U\]\to X$ and
 $s\:Y\to\[V\]$,
 where $U$ and $V$ are simplicial sets.
 Consider the commutative diagram
 $$
 \xymatrix {
 [X,Y]
 \ar[r]^-{h}
 \ar[d]_-{k} &
 [SX,SY]
 \ar[d]^-{l} \\
 [\[U\],\[V\]]
 \ar[r]^-{\~h} &
 [S\[U\],S\[V\]],
 }
 $$
 where
 $\~h$ is the main invariant and
 the map $k$ and the isomorphism $l$ are induced by the pair
 $(r,s)$.
 By Claim~7.2,
 $\~h$ is straight.
 By Lemma~3.1,
 the invariant $\~h\circ k$ is straight.
 Since $h=l^{-1}\circ\~h\circ k$,
 $h$ is also straight.
 \qed
 \end {demo}


 \head {\S~8. The superposition
        $Z\:\<\Si(U,V)\>_0\to\Si(U,\<V\>_0)$}


 For a set $X$,
 let $\<X\>_0\subseteq\<X\>$ be the kernel of the homomorphism
 $\<X\>\to\Z$, $\`x\'\mapsto1$.
 We apply the functor $\<?\>_0$ to simplicial sets degreewise.

 Let $U$ and $V$ be simplicial sets.
 The canonical simplicial map $c=c_V\:V\to\<V\>$ gives rise to
 the map $c^U_\#\:\Si(U,V)\to\Si(U,\<V\>)$ and
 the homomorphism $(c^U_\#)^+\:\<\Si(U,V)\>\to\Si(U,\<V\>)$.
 We have the commutative diagram
 $$
 \xymatrix {
 \<\Si(U,V)\>_0
 \ar[r]^-{Z}
 \ar[d] &
 \Si(U,\<V\>_0)
 \ar[d] \\
 \<\Si(U,V)\>
 \ar[r]^-{(c^U_\#)^+} &
 \Si(U,\<V\>),
 }
 $$
 where
 the vertical arrows are induced by the canonical inclusion
 $\<?\>_0\to\<?\>$ and
 $Z$ is a new homomorphism
 called the {\it superposition}.


 \head {\S~9. Surjectivity of the superposition}


 Our aim here is Lemma~9.1.
 We follow \cite[\S\S\ 12, 13]{me-3}.


 \subhead {Extension of simplicial maps.}
 For $n\ge0$,
 let
 $\sx^n$ be the combinatorial standard $n$-simplex
 (a simplicial set) and
 $\dsx^n$ be its boundary.

 Let $W$ be a contractible fibrant simplicial set.
 For each $n\ge0$,
 choose a map $e_n\:\Si(\dsx^n,W)\to\Si(\sx^n,W)$ such that
 $e_n(q)|_{\dsx^n}=q$
 for any $q\in\Si(\dsx^n,W)$.

 Let $U$ be a simplicial set.
 For each simplicial subset $A\subseteq U$,
 we introduce the map $E_A\:\Si(A,W)\to\Si(U,W)$, $x\mapsto t$,
 where
 $t|_A=x$ and
 $t\circ p=e_n(t\circ p|_{\dsx^n})$
 for the characteristic map $p\:\sx^n\to U$ of each
 non-degenerate simplex outside $A$.
 Clearly,
 \begin {itemize}
 \item [(1)]
 $E_A(x)|_A=x$;
 \item [(2)]
 $E_A(x)|_B=E_{A\cap B}(x|_{A\cap B})|_B$,
 \end {itemize}
 where
 $A,B\subseteq U$ are simplicial subsets and
 $x\in\Si(A,W)$.


 \subhead {The ring $\<Q\>$ and its identity $I$.}
 Let $Q$ be the system of simplicial subsets of $U$
 consisting of
 all subsets isomorphic to $\sx^n$, $n\ge0$, and
 the empty subset.
 Suppose that
 the simplicial set $U$ is
 {\it polyhedral},
 i.~e.\ $Q$ is its cover closed under intersection,
 and
 {\it compact},
 i.~e.\ generated by a finite number of simplices.
 $Q$ is finite.

 We introduce multiplication in $\<Q\>$
 by putting $\`A\'\`B\'=\`A\cap B\'$
 for $A,B\in Q$.
 The ring $\<Q\>$ has an identity $I$.
 Indeed,
 the homomorphism $e\:\<Q\>\to\Z^Q$,
 $$
 e(\`A\')(B)=
 \begin{cases}
 1 &
 \text{if $A\supseteq B$,} \\
 0 &
 \text{otherwise,}
 \end{cases}
 $$
 $A,B\in Q$, is an isomorphism
 (``an upper unitriangular matrix'')
 preserving multiplication.
 Therefore,
 $I=e^{-1}(1)$ is an identity.


 \subhead {The homomorphism
           $K\:\Si(U,\<W\>_0)\to\<\Si(U,W)\>_0$.}
 For a simplicial set $T$,
 let $Z_T\:\<\Si(T,W)\>_0\to\Si(T,\<W\>_0)$ be the
 superposition.
 For simplicial sets $T\supseteq A$,
 let
 $r^T_A\:\Si(T,W)\to\Si(A,W)$ and
 $s^T_A\:\Si(T,\<W\>_0)\to\Si(A,\<W\>_0)$
 be the restriction maps.
 $s^T_A$ is a homomorphism.
 If $T=U$,
 we omit the corresponding sub/superscript in this notation.

 Note that
 $Z_A$ is an isomorphism
 for $A\in Q$.
 Consider the map
 $k\:Q\to\Hom(\Si(U,\<W\>_0),\<\Si(U,W)\>_0)$,
 $A\mapsto\<E_A\>_0\circ Z_A^{-1}\circ s_A$:
 $$
 k(A)\:
 \Si(U,\<W\>_0)
 \xrightarrow{s_A}
 \Si(A,\<W\>_0)
 \xrightarrow{Z_A^{-1}}
 \<\Si(A,W)\>_0
 \xrightarrow{\<E_A\>_0}
 \<\Si(U,W)\>_0.
 $$
 Put $K=k^+(I)$.

 \begin {claim} [9.1. Lemma.]
 The diagram
 $$
 \xymatrix {
 &
 \<\Si(U,W)\>_0
 \ar[d]^-{Z} \\
 \Si(U,\<W\>_0)
 \ar[r]^-{\id}
 \ar[ur]^-{K} &
 \Si(U,\<W\>_0)
 }
 $$
 is commutative.
 \end {claim}

 \begin {demo} [Proof.]
 Take $A,B\in Q$.
 We have the commutative diagram
 $$
 \xymatrix {
 &
 \Si(A,\<W\>_0)
 \ar[r]^{Z_A^{-1}}
 \ar[dd]^-{s^A_C} &
 \<\Si(A,W)\>_0
 \ar[r]^{\<E_A\>_0}
 \ar[dd]^-{\<r^A_C\>_0} &
 \<\Si(U,W)\>_0
 \ar[d]^-{\<r_B\>_0} \\
 \Si(U,\<W\>_0)
 \ar[ur]^-{s_A}
 \ar[dr]_-{s_C} & & &
 \<\Si(B,W)\>_0 \\
 &
 \Si(C,\<W\>_0)
 \ar[r]^{Z_C^{-1}} &
 \<\Si(C,W)\>_0
 \ar[r]^{\<E_C\>_0} &
 \<\Si(U,W)\>_0,
 \ar[u]_-{\<r_B\>_0}
 }
 $$
 where $C=A\cap B$
 (commutativity of the ``pentagon'' follows from the property
 (2) of the family $E$).
 Therefore,
 $\<r_B\>_0\circ k(A)=\<r_B\>_0\circ k(A\cap B)$.
 Therefore,
 $\<r_B\>_0\circ k^+(X)=\<r_B\>_0\circ k^+(X\`B\')$
 for $X\in\<Q\>$.
 We have $
 \<r_B\>_0\circ K=
 \<r_B\>_0\circ k^+(I)=
 \<r_B\>_0\circ k^+(I\`B\')=
 \<r_B\>_0\circ k^+(\`B\')=
 \<r_B\>_0\circ k(B)=
 \<r_B\>_0\circ\<E_B\>_0\circ Z_B^{-1}\circ s_B=
 Z_B^{-1}\circ s_B
 $,
 because $r_B\circ E_B=\id$
 by property (1) of the family $E$.
 We get $
 s_B\circ Z\circ K=
 Z_B\circ\<r_B\>_0\circ K=
 s_B
 $.
 Since $B$ is arbitrary,
 $Z\circ K=\id$.
 \qed
 \end {demo}


 \head {\S~10. A cocartesian square}


 Let
 $U$ be a compact polyhedral simplicial set and
 $V$ be a fibrant simplicial set.
 The canonical simplicial map $c=c_V\:V\to\<V\>$ induces the
 maps
 $c^U_\#\:\Si(U,V)\to\Si(U,\<V\>)$ and
 $c^U_*\:[U,V]\to[U,\<V\>]$.
 Consider the commutative square of abelian groups and
 homomorphisms
 $$
 \xymatrix {
 \<\Si(U,V)\>
 \ar[r]^-{(c^U_\#)^+}
 \ar[d]_-{\<p\>} &
 \Si(U,\<V\>)
 \ar[d]^-{q} \\
 \<[U,V]\>
 \ar[r]^-{(c^U_*)^+} &
 [U,\<V\>],
 }
 $$
 where
 $p=[?]\:\Si(U,V)\to[U,V]$ and
 $q=[?]$
 (the projections).

 \begin {claim} [10.1. Lemma.]
 This square is cocartesian.
 \end {claim}

 \begin {demo} [Proof.]
 Since $\<p\>$ and $q$ are epimorphisms,
 it suffices to show that
 $\Ker q=(c^U_\#)^+(\Ker\<p\>)$.

 Suppose we have a decomposition
 $$
 V=\coprod_{i\in I}V_i.
 $$
 Consider the commutative diagram
 $$
 \xymatrix {
 \bigoplus_{i\in I}\<\Si(U,V_i)\>
 \ar[rrr]^-{\bigoplus_{i\in I}((c_i)^U_\#)^+}
 \ar[ddd]_-{\bigoplus_{i\in I}\<p_i\>}
 \ar[dr] & & &
 \bigoplus_{i\in I}\Si(U,\<V_i\>)
 \ar[ddd]^-{\bigoplus_{i\in I}q_i}
 \ar[dl]^-{E} \\
 &
 \<\Si(U,V)\>
 \ar[r]^-{(c^U_\#)^+}
 \ar[d]_-{\<p\>} &
 \Si(U,\<V\>)
 \ar[d]^-{q} \\
 &
 \<[U,V]\>
 \ar[r]^-{(c^U_*)^+} &
 [U,\<V\>] &
 \\
 \bigoplus_{i\in I}\<[U,V_i]\>
 \ar[rrr]^-{\bigoplus_{i\in I}((c_i)^U_*)^+}
 \ar[ur] & & &
 \bigoplus_{i\in I}[U,\<V_i\>],
 \ar[ul]_-{e}
 }
 $$
 where
 $c_i$, $p_i$, and $q_i$ are similar to $c$, $p$, and $q$
 (respectively) and
 the slanting arrows are induced by the inclusions $V_i\to V$.
 Since $U$ is compact,
 $E$ and $e$ are isomorphisms.
 Therefore,
 is suffices to show that
 $\Ker q_i=((c_i)^U_\#)^+(\Ker\<p_i\>)$
 for each $i\in I$.
 This reduction allows us to assume that
 $V$ is $0$-connected.

 Consider the commutative diagram
 $$
 \xymatrix {
 \<\Si(U,V)\>_0
 \ar[rrr]^-{Z}
 \ar[ddd]_-{\<p\>_0}
 \ar[dr]_-{I} & & &
 \Si(U,\<V\>_0)
 \ar[ddd]^-{q_0}
 \ar[dl]^-{j^U_\#} \\
 &
 \<\Si(U,V)\>
 \ar[r]^-{(c^U_\#)^+}
 \ar[d]_-{\<p\>} &
 \Si(U,\<V\>)
 \ar[d]^-{q} \\
 &
 \<[U,V]\>
 \ar[r]^-{(c^U_*)^+} &
 [U,\<V\>] &
 \\
 \<[U,V]\>_0
 \ar[rrr]^-{z}
 \ar[ur]^-{i} & & &
 [U,\<V\>_0],
 \ar[ul]_-{j^U_*}
 }
 $$
 where
 $q_0=[?]$ (the projection),
 $Z$ is the superposition,
 $z$ is the homomorphism such that
 the outer square is commutative,
 $I$ and $i$ are the inclusion homomorphisms, and
 $j\:\<V\>_0\to\<V\>$ is the inclusion simplicial homomorphism.
 Clearly,
 $\Ker q=j^U_\#(\Ker q_0)$.
 Therefore,
 it suffices to show that
 $\Ker q_0=Z(\Ker\<p\>_0)$.

 Since $V$ is fibrant and $0$-connected,
 there is a surjective simplicial map $f\:W\to V$,
 where $W$ is a contractible fibrant simplicial set.
 Consider the commutative diagram
 $$
 \xymatrix {
 \<\Si(U,W)\>_0
 \ar[r]^-{\~Z}
 \ar[d]_-{\<f^U_\#\>_0} &
 \Si(U,\<W\>_0)
 \ar[d]^-{(\<f\>_0)^U_\#} \\
 \<\Si(U,V)\>_0
 \ar[r]^-{Z}
 \ar[d]_-{\<p\>_0} &
 \Si(U,\<V\>_0)
 \ar[d]^-{q_0} \\
 \<[U,V]\>_0
 \ar[r]^-{z} &
 [U,\<V\>_0],
 }
 $$
 where
 the map $f^U_\#\:\Si(U,W)\to\Si(U,V)$ and
 the simplicial homomorphism $\<f\>_0\:\<W\>_0\to\<V\>_0$
 are induced by $f$ and
 $\~Z$ is the superposition.
 Since $\<f\>_0$ is surjective,
 it is a fibration.
 Therefore,
 $\Ker q_0\subseteq\Im(\<f\>_0)^U_\#$.
 By Lemma~9.1,
 $\~Z$ is surjective.
 Since $W$ is contractible,
 $\Im\<f^U_\#\>_0\subseteq\Ker\<p\>_0$.
 Therefore,
 $\Ker q_0\subseteq Z(\Ker\<p\>_0)$.
 The reverse inclusion is obvious.
 \qed
 \end {demo}


 \head {\S~11. The homomorphism
        $P:\Si(U,\<V\>)\to L(\[U\],\[V\])$}


 For $n\ge0$,
 let
 $\Sx^n$ be the geometric standard $n$-simplex and
 $\iSx^n$ be its interior.
 For
 a simplicial set $U$ and
 a point $z\in\Sx^n$,
 there is a canonical map $z_U\:U_n\to\[U\]$.
 The map $\Sx^n\times U_n\to\[U\]$, $(z,u)\mapsto z_U(u)$, is
 the canonical pairing of geometric realization.

 Let $U$ and $V$ be simplicial sets.
 We define a homomorphism
 $\~P:\Si(U,\<V\>)\to\Hom(\<\[U\]\>,\<\[V\]\>)$.
 For
 $t\in\Si(U,\<V\>)$ and
 $x\in\[U\]$, $x=z_U(u)$,
 where
 $z\in\Sx^n$ and
 $u\in U_n$
 ($n\ge0$),
 put $\~P(t)(\`x\')=\<z_V\>(t_n(u))$:
 $$
 u\in U_n
 \xrightarrow{t_n}
 \<V\>_n=\<V_n\>
 \xrightarrow{\<z_V\>}
 \<\[V\]\>.
 $$
 $\~P$ is well-defined.

 Suppose that
 $U$ is compact.

 \begin {claim} [11.1. Lemma.]
 $\Im\~P\subseteq L(\[U\],\[V\])$.
 \end {claim}

 \begin {demo} [Proof.]
 Let $U^\%_n\subseteq U_n$ ($n\ge0$) be the set of
 non-degenerate simplices.
 For $u\in U^\%_n$ ($n\ge0$),
 we define a homomorphism $I_u\:\<V_n\>\to L(\[U\],\[V\])$.
 For $v\in V_n$, $x\in\[U\]$,
 put
 $$
 I_u(\`v\')(\`x\')=
 \begin{cases}
 \`z_V(v)\' &
 \text{if $x=z_U(u)$ for $z\in\iSx^n$,} \\
 0 &
 \text{otherwise.}
 \end{cases}
 $$
 This equality is preserved if
 we replace
 $\`v\'$ by $w\in\<V_n\>$ and
 $\`z_V(v)\'$ by $\<z_V\>(w)$.
 It suffices to show that
 $$
 \~P(t)=
 \sum_{n\ge0,\,u\in U^\%_n}I_u(t_n(u)),
 \qquad
 t\in\Si(U,\<V\>).
 $$
 Evaluating each side at $\`x\'$, $x=z_U(u)$,
 where
 $z\in\iSx^n$ and
 $u\in U^\%_n$ ($n\ge0$),
 we get $\<z_V\>(t_n(u))$.
 \qed
 \end {demo}

 Lemma~11.1 allows us to introduce the homomorphism
 $P:\Si(U,\<V\>)\to L(\[U\],[V\])$, $P(t)=\~P(t)$.

 \begin {claim} [11.2. Lemma.]
 The diagram
 $$
 \xymatrix {
 \Si(U,V)
 \ar[r]^{c^U_\#}
 \ar[d]_-{\[?\]} &
 \Si(U,\<V\>)
 \ar[d]^-{P} \\
 C(\[U\],\[V\])
 \ar[r]^-{\<?\>} &
 L(\[U\],\[V\]),
 }
 $$
 where
 $c=c_V\:V\to\<V\>$ is the canonical simplicial map,
 is commutative.
 \end {claim}

 \begin {demo} [Proof.]
 For
 $s\in\Si(U,V)$ and
 $x\in\[U\]$, $x=z_U(u)$,
 where
 $z\in\Sx^n$ and
 $u\in U_n$ ($n\ge0$),
 we have
 $
 (P\circ c^U_\#)(s)(\`x\')=
 P(c\circ s)(\`x\')=
 \<z_V\>((c\circ s)_n(u))=
 \`z_V(s_n(u))\'=
 \`\[s\](z_U(u))\'=
 \`\[s\](x)\'=
 \<\[s\]\>(\`x\')
 $.
 \qed
 \end {demo}


 \head {\S~12. Necessity in Theorem~1.1}


 \begin {claim} [12.1. Claim.]
 Let
 $U$ be a compact polyhedral simplicial set,
 $V$ be a fibrant simplicial set,
 $h\:[\[U\],\[V\]]\to[S\[U\],S\[V\]]$ be the main invariant,
 $M$ be an abelian group, and
 $f\:[\[U\],\[V\]]\to M$ be a straight invariant.
 Then there exists
 a homomorphism $d\:[S\[U\],S\[V\]]\to M$ such that
 $f=d\circ h$.
 \end {claim}

 \begin {demo} [Proof.]
 Since $f$ is straight,
 there is a homomorphism $F\:L(\[U\],\[V\])\to M$ such that
 $f([a])=F(\<a\>)$
 for $a\in C(\[U\],\[V\])$.
 Consider the diagram of abelian groups and homomorphisms
 $$
 \xymatrix {
 \<C(\[U\],\[V\])\>
 \ar[rrr]^-{k^+}
 \ar[ddd]_-{\<r\>} & & &
 L(\[U\],\[V\])
 \ar[ddd]^-{F} \\
 &
 \<\Si(U,V)\>
 \ar[r]^-{(c^U_\#)^+}
 \ar[d]_-{\<p\>}
 \ar[ul]^{\<I\>} &
 \Si(U,\<V\>)
 \ar[d]^-{q}
 \ar[ur]_-{P} \\
 &
 \<[U,V]\>
 \ar[r]^-{(c^U_*)^+}
 \ar[dl]_-{\<i\>} &
 [U,\<V\>]
 \ar@{-->}[dr]^-{\~d} &
 \\
 \<[\[U\],\[V\]]\>
 \ar[rrr]^-{f^+} & & &
 M.
 }
 $$
 Here
 the inner square is as in \S~10,
 $r=[?]\:C(\[U\],\[V\])\to[\[U\],\[V\]]$ (the projection),
 $k=\<?\>\:C(\[U\],\[V\])\to L(\[U\],\[V\])$,
 $I=\[?\]\:\Si(U,V)\to C(\[U\],\[V\])$ (the geometric
 realization map),
 $i\:[U,V]\to[\[U\],\[V\]]$, $[s]\mapsto[\[s\]]$, and
 $P$ is as in \S~11.
 By Lemma~11.2,
 the upper trapezium is commutative.
 The solid arrows
 are defined and
 form a commutative subdiagram.
 Since the inner square is cocartesian
 by Lemma~10.1,
 the dashed arrow $\~d$ is well-defined by the condition
 of commutativity of the diagram.

 Consider the diagram
 $$
 \xymatrix {
 \<[U,V]\>
 \ar[rr]^-{(c^U_*)^+}
 \ar[dd]_-{\<i\>} & &
 [U,\<V\>]
 \ar[dd]^-{e}
 \ar[dl]_-{\~d} \\
 &
 M &
 \\
 \<[\[U\],\[V\]]\>
 \ar[rr]^-{h^+}
 \ar[ur]^-{f^+} & &
 [S\[U\],S\[V\]],
 \ar[ul]_-{d}
 }
 $$
 where
 $e$ is the isomorphism from Lemma~4.1 and
 $d=\~d\circ e^{-1}$.
 The square is commutative
 by Lemma~4.1.
 We have $\~d\circ(c^U_*)^+=f^+\circ\<i\>$.
 Since $V$ is fibrant,
 $i$ is a bijection, and
 thus
 $\<i\>$ is an isomorphism.
 We get $f^+=d\circ h^+$
 (so the diagram is commutative).
 Therefore,
 $f=d\circ h$.
 \qed
 \end {demo}

 \begin {claim} [12.2. Proposition.]
 Let
 $X$ be finitary valid space,
 $Y$ be a space,
 $h\:[X,Y]\to[SX,SY]$ be the main invariant,
 $M$ be an abelian group, and
 $f\:[X,Y]\to M$ be a straight invariant.
 Then there exists a homomorphism $d\:[SX,SY]\to M$ such that
 $f=d\circ h$.
 \end {claim}

 \begin {demo} [Proof.]
 There are
 a homotopy equivalence $r\:X\to\[U\]$ and
 a weak homotopy equivalence $s\:\[V\]\to Y$,
 where
 $U$ is a compact polyhedral simplicial set and
 $V$ is a fibrant simplicial set.
 We construct the commutative diagram
 $$
 \xymatrix {
 [\[U\],\[V\]]
 \ar[rr]^-{\~h}
 \ar[dd]_-{k}
 \ar[dr]_-{\~f} & &
 [S\[U\],S\[V\]]
 \ar[dd]^-{l}
 \ar[dl]^-{\~d} \\
 &
 M &
 \\
 [X,Y]
 \ar[rr]^-{h}
 \ar[ur]^-{f} & &
 [SX,SY].
 \ar[ul]_-{d}
 }
 $$
 Here
 the bijection $k$ and
 the isomorphism $l$
 are induced by the pair $(r,s)$ and
 $\~h$ is the main invariant.
 The square is commutative.
 By Lemma~3.1,
 the invariant $\~f=f\circ k$ is straight.
 By Claim~12.1,
 there is a homomorphism $\~d$ such that
 $\~f=\~d\circ\~h$.
 Set $d=\~d\circ l^{-1}$.
 Since $k$ is a bijection,
 we get $f=d\circ h$
 (so the diagram is commutative).
 \qed
 \end {demo}


 \head {\S~13. Three counterexamples}


 \subhead {The Hawaiian ear-ring.}
 Let us show that
 the hypothesis of validity of $Y$ in Theorem~1.1 and
 Proposition~7.3 is essential.
 Let
 $X$ be the one-point compactification of the ray
 $\R_+=(0,\infty)$
 (a circle) and
 $Y$ be that of the space $\R_+\setminus\N$
 (the Hawaiian ear-ring \cite[Example~1.25]{Hatcher}).
 We define a map $m\in C(X,Y)$
 by putting
 $$
 m(x)=
 \left[\frac{x+1}2\right]+
 (-1)^{[x/2]}\{-x\}
 $$
 for $x\in\R_+\setminus\N$.
 Here $[t]$ and $\{t\}$ are the integral and the fractional
 (respectively) parts of a number $t\in\R$.
 The element of $\pi_1(Y,\infty)$ represented by the loop $m$
 is the (reasonably understood) infinite product of commutators
 \begin {equation} \label {*} \tag {\*}
 \prod_{p=0}^\infty[u_{2p},u_{2p+1}],
 \end {equation}
 where $u_q$ is the element realized by the closure of the
 interval $(q,q+1)$.
 Let $e\in H_1(X)$ be the standard generator.
 As in \cite[p.\ 76]{Higman},
 we get that
 the element $m_*(e)\in H_1(Y)$ has infinite order.
 Therefore,
 there is a homomorphism $k\:H_1(Y)\to\Q$ such that
 $k(m_*(e))=1$.
 We define a homomorphism $d\:[SX,SY]\to\Q$
 by putting $d([v])=k(v_*(e))$
 for a morphism $v\:SX\to SY$.
 Let $h\:[X,Y]\to[SX,SY]$ be the main invariant.
 We show that
 {\it the invariants
 $d\circ h$ and
 thus $h$
 are not straight}.

 For $y\in Y$ and $i=0,1$,
 put $y_{(i)}\in Y$ equal
 to $\infty$ if $i=1$ and
 to $y$ otherwise.
 For $i,j=0,1$,
 we define a map $r_{ij}\in C(Y,Y)$.
 For $y\in\R_+\setminus\N$,
 we put $r_{ij}(y)$ equal
 to $y_{(j)}$ if $[y]$ is odd and
 to $y_{(i)}$ otherwise.
 For elements $z_{ij}$, $i,j=0,1$, of an abelian group,
 put $\V_{ij}z_{ij}=z_{00}-z_{10}-z_{01}+z_{11}$.
 Clearly,
 $\V_{ij}\<r_{ij}\>=0$ in $L(Y,Y)$.
 Put $a_{ij}=r_{ij}\circ m\in C(X,Y)$.
 We get $\V_{ij}\<a_{ij}\>=0$ in $L(X,Y)$.
 Therefore,
 $\V_{ij}f([a_{ij}])=0$
 for any straight invariant $f$.
 We show that
 this is false for the invariant $d\circ h$.
 We have $a_{00}=m$;
 the map $a_{11}$ is constant.
 It is easy to see that
 the maps $a_{10}$ and $a_{01}$ are null-homotopic
 (this ``follows formally'' from
 the presentation \eqref{*} and
 the equalities
 $r_{10\,*}(u_{2p})=r_{01\,*}(u_{2p+1})=1$).
 We get $
 \V_{ij}(d\circ h)([a_{ij}])=
 (d\circ h)([m])=
 k(m_*(e))=
 1
 $.
 \qed

 Using \cite[Theorem~2]{Barratt-Milnor},
 one can make the spaces $X$ and $Y$ simply-connected in this
 example.


 \subhead {The Warsaw circle.}
 Let us show that
 the hypothesis of validity of $X$ in Theorem~1.1 and
 Proposition~12.2 is essential.
 Let $X$ be the Warsaw circle
 \cite[Exercise~7 in \S~1.3]{Hatcher} and
 $Y$ be the unit circle in $\C$.
 $Y$ is a topological abelian group.
 The group $[X,Y]$ is
 non-zero
 by \cite[Exercise~7 in \S~1.3, Proposition~1.30]{Hatcher} and
 torsion-free
 by \cite[Theorem~1 in \S~56-III]{Kuratowski}.
 Therefore,
 there is a non-zero homomorphism $f\:[X,Y]\to\Q$.
 By Lemma~7.1,
 $f$ is a straight invariant.
 Since
 $X$ is weakly homotopy equivalent to a point
 \cite[Exercise~10 in \S~4.1]{Hatcher} and
 $Y$ is $0$-connected,
 the main invariant $h\:[X,Y]\to[SX,SY]$ is constant.
 Therefore
 {\it there exists no homomorphism $d\:[SX,SY]\to\Q$ such that
 $f=d\circ h$}.
 \qed


 \subhead {An infinite discrete space.}
 Let us show that the hypothesis of finitarity of $X$ in
 Theorem~1.1 and Proposition~12.2 is essential
 (see also \S~14).

 Note that,
 for an infinite set $X$,
 the subgroup $B(X)\subseteq\Z^X$ is not a direct summand
 because
 the group $\Z^X$ is reduced and
 the group $\Z^X/B(X)$ is divisible and non-zero.

 Let $X$ and $Y$ be discrete spaces,
 $X$ infinite and
 $Y=\{y_0,y_1\}$.
 Introduce the function $k\:Y\to\Z$, $y_i\mapsto i$, $i=0,1$.
 Consider the invariant $f\:[X,Y]\to B(X)$,
 $[a]\mapsto k\circ a$, $a\in C(X,Y)$.

 The invariant $f$ is straight
 because,
 for the homomorphism $F\:L(X,Y)\to B(X)$,
 $F(u)(x)=k^+(u(\`x\'))$, $x\in X$, $u\in L(X,Y)$,
 we have $f([a])=F(\<a\>)$,
 $a\in C(X,Y)$.

 Let $h\:[X,Y]\to[SX,SY]$ be the main invariant.
 We show that
 {\it there exists no homomorphism $d\:[SX,SY]\to B(X)$ such
 that
 $f=d\circ h$}.
 Assume that
 there is such a $d$.

 Consider the homomorphism $l\:\Z^X\to\Hom(\<X\>,\<Y\>)$,
 $l(v)(\`x\')=v(x)(\`y_1\'-\`y_0\')$, $x\in X$, $v\in\Z^X$.
 We have $l(f([a]))=\<a\>-\<a_0\>$,
 $a\in C(X,Y)$,
 where $a_0\:X\to Y$, $x\mapsto y_0$.
 Clearly,
 there is an isomorphism $e\:\Hom(\<X\>,\<Y\>)\to[SX,SY]$
 such that
 $e(\<a\>)=h([a])$,
 $a\in C(X,Y)$.
 Consider the composition
 $$
 r\:
 \Z^X
 \xrightarrow{l}
 \Hom(\<X\>,\<Y\>)
 \xrightarrow{e}
 [SX,SY]
 \xrightarrow{d}
 B(X).
 $$
 For $a\in C(X,Y)$,
 we have $
 r(f([a]))=
 (d\circ e\circ l\circ f)([a])=
 d(e(\<a\>-\<a_0\>))=
 d(h([a])-h([a_0]))=
 f([a])-f([a_0])=
 f([a])
 $.
 Since the elements $f([a])$, $a\in C(X,Y)$, generate $B(X)$,
 we get $r|_{B(X)}=\id$,
 which is impossible.
 \qed


 \head {\S~14. Invariants of maps $\R P^\infty\to\R P^\infty$}


 Here we show that
 the hypothesis of finitarity of $X$ in Theorem~1.1 and
 Proposition~12.2 is essential even if $M$ is divisible.
 (Possibly,
 if
 $M$ is divisible and/or
 $Y$ is (simply-)connected,
 the hypothesis of finitarity of $X$ can be replaced by the
 weaker one that
 $Y$ is weakly homotopy equivalent to a finite-dimensional
 CW-complex.)

 Let $X$ and $Y$ be spaces.
 A set $E\subseteq X$ is called {\it $Y$-representative\/} if
 any maps $a,b\in C(X,Y)$ equal on $E$ are homotopic.
 $X$ is called {\it $Y$-unitary\/} if
 any finite cover of $X$ contains a $Y$-representative set.

 \begin {claim} [14.1. Lemma.]
 Let $M$ be a divisible group.
 If $X$ is $Y$-unitary,
 then any invariant $f\:[X,Y]\to M$ is straight.
 \end {claim}

 \begin {demo} [Proof.]
 Introduce the maps
 $r=[?]\:C(X,Y)\to[X,Y]$
 (the projection) and
 $k=\<?\>\:C(X,Y)\to L(X,Y)$.
 We seek a homomorphism $F$ giving the commutative diagram
 $$
 \xymatrix {
 \<C(X,Y)\>
 \ar[r]^-{k^+}
 \ar[d]_-{\<r\>} &
 L(X,Y)
 \ar[d]^-{F} \\
 \<[X,Y]\>
 \ar[r]^-{f^+} &
 M.
 }
 $$
 Since $M$ is divisible,
 it suffices to show that
 $\Ker k^+\subseteq\Ker\<r\>$.
 Take an element $w\in\Ker k^+$.
 We show that
 $w\in\Ker\<r\>$.
 There are
 a finite set $I$,
 a map $l\:I\to C(X,Y)$, and
 an element $v\in\<I\>$
 such that
 $\<l\>(v)=w$.
 Put $a_i=l(i)$, $i\in I$.
 For an equivalence $d$ on $I$,
 let $p_d\:I\to I/d$ be the projection.
 Let $N$ be the set of equivalences $d$ on $I$ such that
 $\<p_d\>(v)=0$ in $\<I/d\>$.

 Take $x\in X$.
 Consider the equivalence $d(x)=\{(i,j):a_i(x)=a_j(x)\}$ on
 $I$.
 We show that
 $d(x)\in N$.
 We have the commutative diagrams
 $$
 \xymatrix {
 I
 \ar[r]^-{l}
 \ar[d]_-{p_{d(x)}} &
 C(X,Y)
 \ar[d]^-{e_x} &
 \<C(X,Y)\>
 \ar[r]^-{k^+}
 \ar[d]_-{\<e_x\>} &
 L(X,Y)
 \ar[dl]^-{h_x} \\
 I/d(x)
 \ar[r]^-{l_x} &
 Y, &
 \<Y\>, &
 }
 $$
 where
 the map $l_x$ is defined by the condition of commutativity of
 the diagram,
 $e_x$ is the map of evaluation at $x$, and
 $h_x$ is the homomorphism of evaluation at $\`x\'$.
 We get $
 \<l_x\>(\<p_{d(x)}\>(v))=
 \<e_x\>(\<l\>(v))=
 \<e_x\>(w)=
 h_x(k^+(w))=
 0
 $.
 Since $l_x$ is injective,
 we get $\<p_{d(x)}\>(v)=0$,
 which is what we promised.

 For an equivalence $d$ on $I$,
 put $E_d=\{x\in X:(i,j)\in d\Rightarrow a_i(x)=a_j(x)\}$.
 Since $x\in E_{d(x)}$
 for any $x\in X$,
 the family $E_d$, $d\in N$, is a cover of $X$.
 Since $X$ is $Y$-unitary,
 $E_d$ is $Y$-representative
 for some $d\in N$.
 For $(i,j)\in d$,
 the maps $a_i$ and $a_j$ are
 equal on $E_d$ and
 thus homotopic.
 Therefore,
 there is a map $m$ giving the commutative diagram
 $$
 \xymatrix {
 I
 \ar[r]^-{l}
 \ar[d]_-{p_d} &
 C(X,Y)
 \ar[d]^-{r} \\
 I/d
 \ar[r]^-{m} &
 [X,Y].
 }
 $$
 We get $
 \<r\>(w)=
 \<r\>(\<l\>(v))=
 \<m\>(\<p_d\>(v))=
 0
 $
 because $d\in N$.
 \qed
 \end {demo}

 Hereafter,
 let $X$ and $Y$ be homeomorphic to $\R P^\infty$.

 \begin {claim} [14.2. Lemma.]
 $X$ is $Y$-unitary.
 \end {claim}

 \begin {demo} [Proof.]
 Let $H^\bullet$ be the $\Z_2$-cohomology.
 Let $g\in H^1X$ and $h\in H^1Y$ be the non-zero classes.

 We show that
 \0
 a set $E\subseteq X$ is $Y$-representative if
 $g|_U\ne0$
 for any neighbourhood $U$ of $E$.
 If maps $a,b\in C(X,Y)$ are equal on $E$,
 they are homotopic on some neighbourhood $U$ of $E$.
 Then
 $a^*(h)|_U=b^*(h)|_U$.
 Since $g|_U\ne0$,
 the homomorphism $?|_U\:H^1X\to H^1U$ is injective.
 Therefore,
 $a^*(h)=b^*(h)$.
 Since $Y$ is a $\mathcal K(\Z_2,1)$ space,
 $a$ and $b$ are homotopic,
 as needed.

 We show that
 $X$ is $Y$-unitary.
 Assume that
 $X=E_1\cup\dotsc\cup E_n$,
 where the sets $E_i$ are not $Y$-representative.
 By \0,
 each $E_i$ has a neighbourhood $U_i$ with $g|_{U_i}=0$.
 Since $U_1\cup\dotsc\cup U_n=X$,
 we get $g^n=0$,
 which is false.
 \qed
 \end {demo}

 We have $[X,Y]=\{u_0,u_1\}$,
 where
 $u_0$ is the class of a constant map and
 $u_1$ is that of a homeomorphism.
 Consider the invariant $f\:[X,Y]\to\Q$, $u_i\mapsto i$,
 $i=0,1$.
 By Lemmas 14.2 and 14.1,
 $f$ is straight.
 Let $h\:[X,Y]\to[SX,SY]$ be the main invariant.
 Using the isomorphism
 $$
 [SX,SY]
 \longrightarrow
 \prod_{i\in\Z}\Hom(H_iX,H_iY),
 \qquad
 [v]\mapsto v_*,
 $$
 we get $2h(u_0)=2h(u_1)$.
 Therefore,
 {\it there exists no homomorphism $d\:[SX,SY]\to\Q$ such that
 $f=d\circ h$}.
 \qed


 \head {\S~15. $K$-straight invariants}


 Let $K$ be a unital ring.
 {\it $K$-modules\/} are unital.


 \subhead {$K$-module $L_K(X,Y)$.}
 For a set $X$,
 let $\<X\>_K$ be the (free) $K$-module with the basis
 $X^\=_K\subseteq\<X\>_K$
 endowed with the bijection $X\to X^\=_K$, $x\mapsto\`x\'_K$.
 For sets $X$ and $Y$,
 let $L_K(X,Y)\subseteq\Hom_K(\<X\>_K,\<Y\>_K)$ be the
 $K$-submodule generated by the $K$-homomorphisms $u$ such that
 $u(X^\=_K)\subseteq Y^\=_K\cup\{0\}$.
 A map $a\:X\to Y$ induces a $K$-homomorphism
 $\<a\>_K\in L_K(X,Y)$, $\<a\>_K(\`x\'_K)=\`a(x)\'_K$.


 \subhead {$K$-straight invariants.}
 Let
 $X$ and $Y$ be spaces and
 $M$ be a  $K$-module.
 An invariant $f\:[X,Y]\to M$ is called {\it $K$-straight\/} if
 there exists a $K$-homomorphism $\~F\:L_K(X,Y)\to M$ such that
 $f([a])=\~F(\<a\>_K)$
 for all $a\in C(X,Y)$.

 \begin {claim} [15.1. Proposition.]
 An invariant $f\:[X,Y]\to M$ is $K$-straight
 if and only if
 it is straight.
 \end {claim}

 Proof is given in \S~16.


 \subhead {The $K$-main invariant
           $\~h\:[X,Y]\to[S_KX,S_KY]_K$.}
 Let
 $S_KX$ be the $K$-complex of singular chains of $X$ with
 coefficients in $K$ and
 $[S_KX,S_KY]_K$ be the $K$-module of $K$-homotopy classes of
 $K$-morphisms $S_KX\to S_KY$.
 For $a\in C(X,Y)$,
 let
 $S_Ka\:S_KX\to S_KY$ be the induced $K$-morphism and
 $[S_Ka]_K\in[S_KX,S_KY]_K$ be its $K$-homotopy class.
 The invariant $\~h\:[X,Y]\to[S_KX,S_KY]_K$,
 $[a]\mapsto[S_Ka]_K$, is called {\it $K$-main}.

 \begin {claim} [15.2. Theorem.]
 Suppose that
 $X$ is valid and finitary and
 $Y$ is valid.
 An invariant $f\:[X,Y]\to M$ is $K$-straight
 if and only if
 there exists a $K$-homomorphism $\~d\:[S_KX,S_KY]_K\to M$ such
 that
 $f=\~d\circ\~h$.
 \end {claim}

 Proof is given in \S~16.
 For $K=\Z$,
 this is Theorem~1.1.


 \head {\S~16. $K$-straight invariants: proofs}


 Let $X$ and $Y$ be sets.
 We define a homomorphism $e\:L(X,Y)\to L_K(X,Y)$.
 For $u\in L(X,Y)$,
 let $e(u)$ be the $K$-homomorphism giving the commutative
 diagram
 $$
 \xymatrix {
 \<X\>
 \ar[r]^-{u}
 \ar[d]_-{i_X} &
 \<Y\>
 \ar[d]^-{i_Y} \\
 \<X\>_K
 \ar[r]^-{e(u)} &
 \<Y\>_K,
 }
 $$
 where
 $i_X$ is the homomorphism $\`x\'\mapsto\`x\'_K$ and
 $i_Y$ is similar.

 For
 an abelian group $A$,
 a $K$-module $M$, and
 a homomorphism $t\:A\to M$,
 we introduce the $K$-homomorphism $t^{(K)}\:K\otimes A\to M$,
 $1\otimes a\mapsto t(a)$.

 \begin {claim} [16.1. Lemma.]
 $e^{(K)}\:K\otimes L(X,Y)\to L_K(X,Y)$ is a $K$-isomorphism.
 \end {claim}

 \begin {demo} [Proof.]
 For $w\in\<Y\>_K$ and $y\in Y$,
 let $w/y\in K$ be the coefficient of $\`y\'_K$ in $w$.
 For $v\in L_K(X,Y)$ and $k\in K\setminus\{0\}$,
 we introduce the homomorphism $v_k\in L(X,Y)$,
 $$
 v_k(\`x\')=
 \sum_{y\in Y:v(\lq x\rq_K)/y=k}\`y\',
 \qquad
 x\in X.
 $$
 It is not difficult to verify that
 the map $d\:L_K(X,Y)\to K\otimes L(X,Y)$,
 $$
 d(v)=
 \sum_{k\in K\setminus\{0\}}k\otimes v_k,
 $$
 is a $K$-homomorphism.
 Using this,
 we get
 $e^{(K)}\circ d=\id$ and
 $d\circ e^{(K)}=\id$.
 \qed
 \end {demo}


 \subhead {Proof of Proposition~15.1.}
 Necessity.
 Let $f$ be $K$-straight.
 There is a $K$-homomorphism $\~F\:L_K(X,Y)\to M$ such that
 $f([a])=\~F(\<a\>_K)$, $a\in C(X,Y)$.
 Consider the homomorphism $F=\~F\circ e$:
 $$
 \xymatrix {
 C(X,Y)
 \ar[rr]^{\<?\>_K}
 \ar[dd]_-{[?]}
 \ar[dr]_-{\<?\>} & &
 L_K(X,Y)
 \ar[dd]^-{\~F} \\
 &
 L(X,Y)
 \ar[ur]_-{e}
 \ar[dr]^-{F} &
 \\
 [X,Y]
 \ar[rr]^-{f} & &
 M.
 }
 $$
 The diagram is commutative.
 We get $f([a])=F(\<a\>)$, $a\in C(X,Y)$.
 Therefore,
 $f$ is straight.

 Sufficiency.
 Let $f$ be straight.
 There is a homomorphism $F\:L(X,Y)\to M$ such that
 $f([a])=F(\<a\>)$, $a\in C(X,Y)$.
 By Lemma~16.1,
 $e^{(K)}$ is a $K$-isomorphism.
 Consider the homomorphism $\~F=F^{(K)}\circ(e^{(K)})^{-1}$:
 $$
 \xymatrix {
 C(X,Y)
 \ar[rr]^{\<?\>_K}
 \ar[dd]_-{[?]}
 \ar[dr]_-{\<?\>} & &
 L_K(X,Y)
 \ar@/^8ex/[dd]^-{\~F} \\
 &
 L(X,Y)
 \ar[r]^-{1\otimes?}
 \ar[ur]^-{e}
 \ar[dr]_-{F} &
 K\otimes L(X,Y)
 \ar[u]_-{e^{(K)}}
 \ar[d]^-{F^{(K)}} &
 \\
 [X,Y]
 \ar[rr]^-{f} & &
 M.
 }
 $$
 The diagram is commutative.
 We get $f([a])=\~F(\<a\>_K)$,
 $a\in C(X,Y)$.
 Therefore,
 $f$ is $K$-straight.
 \qed


 \subhead {The homomorphism $I\:[SX,SY]\to[S_KX,S_KY]_K$.}
 Let $X$ and $Y$ be spaces.
 A morphism $v\:SX\to SY$ induces a $K$-morphism
 $$
 S_KX=
 K\otimes SX
 \xrightarrow{\id\otimes v}
 K\otimes SY=
 S_KY.
 $$
 Consider the homomorphism $I\:[SX,SY]\to[S_KX,S_KY]_K$,
 $[v]\mapsto[\id\otimes v]_K$.

 \begin {claim} [16.2. Lemma.]
 If the group $H_\bullet(X)$ is finitely generated,
 then the $K$-homomorphism
 $$
 I^{(K)}\:K\otimes[SX,SY]\to[S_KX,S_KY]_K
 $$
 is a $K$-split $K$-monomorphism,
 i.~e.\ there exists a $K$-homomorphism
 $R\:[S_KX,S_KY]_K\to K\otimes[SX,SY]$ such that
 $R\circ I^{(K)}=\id$.
 \end {claim}

 \begin {dem}
 This is a variant of the universal coefficient theorem,
 cf.\ \cite[Theorems 5.2.8 and 5.5.10]{Spanier}.
 \qed
 \end {dem}


 \subhead {Proof of Theorem~15.2.}
 We have $\~h=I\circ h$,
 where $h\:[X,Y]\to[SX,SY]$ is the main invariant.
 By Proposition~7.3,
 $h$ is straight.
 Therefore,
 $\~h$ is straight.
 By Proposition~15.1,
 $\~h$ is $K$-straight.

 This gives the sufficiency.
 Necessity.
 Let $f$ be $K$-straight.
 By Proposition~15.1,
 $f$ is straight.
 By Proposition~12.2,
 there is a homomorphism $d\:[SX,SY]\to M$ such that
 $f=d\circ h$.
 By Lemma~16.2,
 there is a $K$-homomorphism $\~d$ such that
 $\~d\circ I^{(K)}=d^{(K)}$:
 $$
 \xymatrix {
 &
 [SX,SY]
 \ar@/^30ex/[dd]^-{I}
 \ar[d]^-{d}
 \ar[dr]_-{1\otimes?} &
 \\
 [X,Y]
 \ar[r]^-{f}
 \ar[ur]^-{h}
 \ar[dr]_-{\~h} &
 M &
 K\otimes[SX,SY]
 \ar[l]_-{d^{(K)}}
 \ar[dl]_-{I^{(K)}} &
 \\
 &
 [S_KX,S_KY]_K.
 \ar[u]_-{\~d} &
 }
 $$
 The diagram is commutative.
 In particular,
 $f=\~d\circ\~h$.
 \qed



 {\noindent \tt ssp@pdmi.ras.ru}

 {\noindent \tt http://www.pdmi.ras.ru/\${}ssp}

 \end {document}